\documentclass[11pt,a4paper]{amsart}
\usepackage{amssymb}
\usepackage[cp1250]{inputenc}
\usepackage{amsmath}
\usepackage{amsfonts}
\usepackage{amsthm}
\usepackage{bbding}
\usepackage{ifsym}
\usepackage{eufrak}
\usepackage{slashbox}
\usepackage{hyperref}
\usepackage[mathscr]{euscript}
\title[FIRST ORDER SYSTEMS OF ODES]{FIRST ORDER SYSTEMS OF ODES WITH NONLINEAR NONLOCAL BOUNDARY CONDITIONS}
\address{Technical University of \L\'od\'z\\
Institute of Mathematics\\
ul.\ W\'ol\-cza\'n\-ska~215\\
93--005 \L\'od\'z\\
POLAND}
\author{Igor Kossowski and Bogdan Przeradzki}
\email{kossowski.igor@gmail.com}
\email{bogdan.przeradzki@p.lodz.pl}
\keywords{Fredholm operator index zero, Mawhin theorem, nonlocal BVP}
\newtheorem{thm}{Theorem}
\newtheorem{lem}[thm]{Lemma}

\theoremstyle{definition}
\newtheorem{rem}{Remark}

\newcommand{\jnt}{\int\limits}
\newcommand{\pog}{\mathbf}
\newcommand{\bb}{\mathbb}

\newcommand{\kr}{\cdot}

\newcommand{\cl}{\overline}
\newcommand{\Om}{\Omega}
\newcommand\ip[1]{\left\langle#1\right\rangle}
\newcommand\norm[1]{\left\|#1\right\|}

\newcommand{\de}{\partial}

\newcommand{\pal}{\mbox{\boldmath $\alpha$}}

\DeclareMathOperator{\ind}{ind}

\DeclareMathOperator{\im}{im}

\DeclareMathOperator{\sgn}{sgn}

\DeclareMathOperator{\codim}{codim}
\DeclareMathOperator{\dom}{dom}
\DeclareMathOperator{\var}{var}
\begin{document}
\begin{abstract}
In this article, we prove an existence of solutions for a nonlocal boundary value problem with nonlinearity in a nonlocal condition. Our method is based upon the Mawhin's coincidence theory.
\end{abstract}\maketitle
\section{Introduction}
In this paper we consider the following ordinary differential equation 
\begin{equation}
x'=f(t,x) \label{eq:ODE}
\end{equation}
with the non-local condition
\begin{equation}
h\Bigg(\jnt_{0}^{1}x(s)\,dg(s)\Bigg)=0 \label{eq:NC},
\end{equation}
where $f:[0,1]\times\bb{R}^{k}\to\bb{R}^{k}$ is continuous, $g=(g^{1},\ldots,g^{k}):[0,1]\to\bb{R}^{k}$ has bounded variation, $h:\bb{R}^{k}\to\bb{R}^{k}$ is continuous and 
\[\jnt_{0}^{1}x(s)\,dg(s)=\Bigg(\jnt_{0}^{1}x^{1}(s)\,dg^{1}(s),\ldots,\jnt_{0}^{1}x^{k}(s)\,dg^{k}(s)\Bigg). \]\par 
The subject of nonlocal boundary conditions for ordinary differential equations has been a topic of various studies in mathematical articles for many years. The multi-point conditions such as $F(x(t_{1}), x(t_{2}),..., x(t_{n}))=0$ were studied at first (\cite{PS}), then also the significantly nonlocal conditions with the values of the unknown function occurring over the entire domain (integral) became the subject of interest. It is easy to see that the conditions which there is the Stieltjes integral with respect to any function with the total variation in contain also multi-point problems.
\par
Usually the matter of consideration are the second-order differential equations because of their supposed applications but sometimes also the first-order differential equations are being considered as in the present paper (\cite{FNR} and \cite{Z}). And with our level of generality the second order differential equations can be treated as the first-order systems. The methods are typical: searching for the fixed point of integral operator using Contraction Principle, Schauder fixed-point theorem, topological-order methods, e.g. basing on Cone Expansion and Compression Theorem, or finally the Leray-Schauder degree of compact mapping or the Mawhin degree of coincidence. \par 
In this paper both differential equations and boundary conditions are nonlinear what somehow forces to the use of the degree of coincidence  - the linear part $x'$ has the nontrivial kernel. Using this method and with such a generality of assumptions the theorems that can be obtained are the ones in which the Brouwer degree of the nonlinear part being not zero on the kernel of the linear part is the main assumption. In this paper there is only the degree of "the half" of the nonlinear operator, i.e. $h$ and the assumptions regarding the other half of the nonlinear part are different. \par 
Nonlinear boundary conditions have occurred before in works \cite{I1}, \cite{I2}, \cite{YRA}, \cite{Y} but they were of different nature than here: under Stieltjes integral there was the assumption of the unknown function with the nonlinear function. Therefore, the obtained results are not comparable with the previous works; those results present a new direction of research. It is possible only to notice the compatibility with the conventional results regarding the existence of the periodic solutions (\cite{M}). This problem will be explained further in paragraph 4. 
\par Let us present a few problems that are similar though different to (\ref{eq:ODE}), (\ref{eq:NC}).  Our problem includes linear non-local condition - $\jnt_{0}^{1}x(s)\,dg(s)=0$.\ There are many papers investigating BVPs with linear non-local conditions (compare with \cite{AB}, \cite{KT1}, \cite{KT2}, \cite{KT3}, \cite{SD}, \cite{W} and the references therein). Our result includes the BVP in \cite{KT3}, namely 
\[ x''=f(t,x,x'), \quad x(0)=0, \quad \jnt_{0}^{1}x'(s)\,dg(s)=0,\] which is at resonance, then $g(1)-g(0)=0$.
\par  In \cite{I1} and \cite{YRA}, authors considered the existence of positive solutions of nonlinear nonlocal BVP of the form $-x''(t)=q(t)f(t,x(t))$ with integral boundary conditions. G. Infante studied with nonlinear integral boundary conditions (see \cite{I1})
\[ x'(0)+H_{1}\Bigg(\jnt_{0}^{1}x(s)\,dA(s)\Bigg)=0, \quad \sigma x'(1)+x(\eta)=H_{2}\Bigg(\jnt_{0}^{1}x(s)\,dB(s)\Bigg).\]
In \cite{YRA}, authors considered another kind of boundary conditions, namely \[x(0)=\jnt_{0}^{1}(x(s))^{a}\,dA(s), \quad x(1)=\jnt_{0}^{1}(x(s))^{b}\,dB(s),\] where $a,b\geq 0$. 
\section{Some preliminaries}
In this section we recall some facts about a Fredholm operator and Mawhin's coincidence theory. This section is based on \cite{RGJM} (page 10-40). \par 
Let $X$ and $Y$ be a Banach space. An linear operator $L:X\supset \dom{L}\to Y$ is said to be a \emph{Fredholm operator} if $\dim\ker{L}<\infty$, $\im{L}$ is closed in $Y$ and $\codim{\im{L}}<\infty$. The index of the Fredholm operator is defined as follows \[\ind{L}:=\dim\ker{L}-\codim{\im{L}}.\] If $L$ is the Fredholm operator, then continuous projections $P:X\to X$, $Q:Y\to Y$ such that $\im{P}=\ker{L}$, $\ker{Q}=\im{L}$ exist. Thus 
$X=\ker{L}\oplus\ker{P}$ and  $Y=\im{L}\oplus \im{Q}$. It is apparent that $\ker{L}\cap\ker{P}=\{0\}$, therefore  we can consider the restriction $L_{P}:=L|_{\ker{P}}:\dom{L}\cap\ker{P}\to Y$ which is invertible. 
A nonlinear operator $N:X\to Y$ is called \emph{$L$-compact} if N maps bounded sets into bounded ones and $K_{P,Q}=L_{P}^{-1}(I_{Y}-Q)$ (by $I_{Y}$ we denote the identity on $Y$) is completely continuous. \par 
Let $L:X\supset \dom{L}\to Y$ be a Fredholm operator of index zero. Since $\dim\ker{L}=\codim\im{L}$ there exists an isomorphism $J:\im{Q}\to\ker{L}$. \par To obtain the results of the existence we use the following Mawhin's theorem. 
\begin{thm}\textup{(\textbf{Mawhin Continuation Theorem})}
 Let $\Omega$ be a bounded open set in $X$. Assume that 
$L:X\supset \dom{L}\to Y$ is a Fredholm operator with index zero and $N$ is $L$-compact. If: 
\begin{enumerate}
\item  equations $Lx=\lambda N(x)$ have no solutions $x\in \dom{L}\cap\partial\Omega$ for all $\lambda\in(0,1]$;
\item Brouwer degree (\cite{D}, page 1-17) $\deg(JQN,\ker{L}\cap \Omega,0)\neq 0$, which  is called \emph{coincidence degree $L$ and $N$}. .\end{enumerate}
Then the equation $Lx=N(x)$ has a solution in $\overline{\Omega}$.
\end{thm} 
Now we return to main problem and present our notations. 
We set $X:=C([0,1],\bb{R}^{k})$, $\dom{L}:=C^{1}([0,1],\bb{R}^{k})$, $Y:=C([0,1],\bb{R}^{k})\times\bb{R}^{k}$
and define mappings $L:X\supset\dom{L}\to Y$, $N:X\to Y$:
$Lx:=(x',0)$ for $x\in\dom{L}$, $(\forall{t\in[0,1]})\,  N(x)=(F(x),h(\jnt_{0}^{1}x(s)\,dg(s)))$ for $x\in X$, where $F$ is Nemytskii operator, i.e. $(F(x))(t)=f(t,x(t))$ for $t\in[0,1]$.
Thus, we obtain 
\begin{equation}
Lx=N(x).\label{eq:ANE}
\end{equation}
It is clear that
\[\ker{L}=\{x\in C^{1}([0,1],\bb{R}^{k}):x=\textup{const.}\},\]
hence $\dim\ker{L}=k<\infty$. Also observe that $\im{L}=C([0,1],\bb{R}^{k})\times\{0\}$, so $\codim{\im{L}}=k$. Consequently, $L$ is a Fredholm operator of index zero and we can use Mawhin's 
theory.
\section{The existence of solutions}
We know that our operator $L$ is a Fredholm operator with index zero. Our purpose is to use the Mawhin's theory. In first step we define projections $P:X\to X$ by
\[ (\forall t\in[0,1])\,(Px)(t):=x(0) \quad \text{for} \quad x\in X,\] and 
$Q:Y\to Y$ by 
\[ Q(z,\pal):=(-\pal,\pal) \quad \text{for} \quad (z,\pal)\in Y.\]
The description of P makes it evident that $\ker{P}=\{x\in X:x(0)=0\}$, hence 
\[ \dom{L}\cap \ker{L}=\{x\in C^{1}([0,1],\bb{R}^{k}):x(0)=0\}.\] Then inverse operator is defined as 
\[ L^{-1}_{P}(z,0)(t)=\jnt_{0}^{t}z(s)\,ds \quad \text{for} \quad z\in C([0,1],\bb{R}^{k})\]  
and we have 
\[ K_{P,Q}(z,\pal)(t)=L^{-1}_{P}(I-Q)(z,\pal)(t)=\jnt_{0}^{t}z(s)\,ds+t\pal \quad \text{for} \quad (z,\pal)\in Y.\]  
Therefore 
\[ (K_{P,Q}N)(x)(t)=\jnt_{0}^{t}f(s,x(s))\,ds+th\Bigg(\jnt_{0}^{1}x(s)\,dg(s)\Bigg).\]
Since the first term is a composition of  Nemytskii operator and Volterra integral operator and the second term is a finite rank we get the following
\begin{lem}
The operator $K_{P,Q}N:X\to Y$ is completely continuous. Therefore operator $N:X\to Y$ is $L$-compact.
\end{lem}
Our main result is given in the following theorem
\begin{thm}
Let us assume that $g(0^{+})\neq g(0)$ and $\lim\limits_{\varepsilon\to 0^{+}}\var(g,[\varepsilon,1])\leq\min\limits_{j\leq k}|g_{j}(0^{+})-g_{j}(0)|$. Then the BVP (\ref{eq:ODE}), (\ref{eq:NC}) has at least one solution if there exists $R>0$ such that:
\begin{enumerate}
\item[(i)] $\ip{f(t,x),x}\leq 0$ for $t\in(0,1]$, $|x|=R$. 
\item[(ii)] Let \[r_{-}:=R\big(\min_{j\leq k}|g^{j}(0^{+})-g^{j}(0)|-\lim\limits_{\varepsilon\to 0^{+}}\var(g,[\varepsilon,1])\big),\] \[r_{+}:=R\big(|g(0^{+})-g(0)|+\lim\limits_{\varepsilon\to 0^{+}}\var(g,[\varepsilon,1])\big).\] $h(x)\neq 0$ for $r_{-}<|x|\leq r_{+}$ and the Brouwer degree $\deg(h,B_{\bb{R}^{k}}(0,r),0)$ is defined and does not vanish for some $r\in(r_{-},r_{+}]$. 
\end{enumerate}
\label{thm:1thm}
\end{thm}
Before we proceed to the proof we recall some notions regarding the  Riemann-Stieltjes integral (\cite{RN}, page 9-11  and 105-123). Let $g:[a,b]\to\bb{R}^{k}$ and consider the sum \[\sum\limits_{i=1}^{n}|g(s_{i})-g(s_{i-1})|,\]
where $a=s_{0}<\ldots<s_{n}=b$. The supremum taken over the set of all partitions of the interval $[a,b]$ is called \emph{the total variation} of the function $g$ on $[a,b]$, which is denoted by $\var(g,[a,b])$.
\begin{lem}
For any continuous function $\varphi:[0,1]\to\bb{R}^{k}$, 
\[ \jnt_{0}^{1}\varphi(s)\,dg(s)=\varphi(0)(g(0^{+})-g(0))+\lim_{\varepsilon\to 0^{+}}\jnt_{\varepsilon}^{1}\varphi(s)\,dg(s)\] 
and the norm of the integral is bounded 
\[ \Bigg|\jnt_{\varepsilon}^{1}\varphi(s)\,dg(s)\Bigg|\leq \sup_{s\in[\varepsilon,1]}|\varphi(s)|\kr \var(g,[\varepsilon,1]).\]
\end{lem}
Recall that both sides expressions are vectors, which means that the first summand has coordinates 
\[\varphi(0)(g(0^{+})-\varphi(0))=\Big(\varphi^{1}(0)\big(g^{1}(0^{+})-g^{1}(0\big)\big),\ldots, \varphi^{k}(0)\big(g^{k}(0^{+})-g^{k}(0)\big)\Big).\]
The proof follows from the form of Riemmann-Stieltjes sums which converge to the integrals:
\[ \Bigg|\sum_{j=1}^{k}\varphi^{j}(s_{j})(g^{j}(t_{j})-g^{j}(t_{j-1}))\Bigg|\leq \sum_{j=1}^{k}|\varphi^{j}(s_{j})|\kr|g^{j}(t_{j})-g^{j}(t_{j-1})|\]
\[\leq \sup_{s\in[\varepsilon,1]}|\varphi(s)|\kr \sum_{j=1}^{k}|g^{j}(t_{j})-g^{j}(t_{j-1})|\leq \sup_{s\in[\varepsilon,1]}|\varphi(s)|\kr \var(g,[\varepsilon,1]).\]
\begin{proof} The proof is carried out in two steps. In step 1, we prove that BVP (\ref{eq:ODE}, \ref{eq:NC}) has the solution under stronger assumptions: $\lim\limits_{\varepsilon\to 0^{+}}\var(g,[\varepsilon,1])<|g(0^{+})-g(0)|$ and 
$\ip{f(t,x),x}<0$ for $t\in(0,1]$, $|x|=R$. \par
\textbf{Step 1.} We know that the BVP (\ref{eq:ODE}), (\ref{eq:NC}) is equivalent to (\ref{eq:ANE}). A linear operator $L$ is a Fredholm operator with index zero and nonlinear $N$ is $L$-compact. If we prove other assumptions of Mawhin theorem  we get the assertion. 
\par Let us consider the family of equations $Lx=\lambda N(x)$, where $\lambda\in(0,1]$. Thus we have the family of problems 
\begin{equation}
\left\{\begin{array}{l}
x'=\lambda f(t,x), \\ 
h\Big(\jnt_{0}^{1}x(s)\,dg(s)\Big)=0.
\end{array}\right. \label{eq:FBVP}
\end{equation}
Now, we shall show that BVPs (\ref{eq:FBVP}) have no solution in $\de\Omega=\de B_{C([0,1],\bb{R}^{k})}(0,R)$ for $\lambda\in(0,1]$. Let us suppose that there exists a solution $\varphi$ of the (\ref{eq:FBVP}) such that $\norm{\varphi}_{C([0,1],\bb{R}^{k})}=R$. We consider then a function
$\psi(t):=|\varphi(t)|^{2}$. Let us assume that $\psi(t_{0})=R^{2}$ for some $t_{0}\in (0,1]$. Then, by the assumption (i), since $\varphi$ is a solution of (\ref{eq:FBVP}) and $|\varphi(t_{0})|=R$, we get a contradiction. Indeed, we obtain 
\[ 0\leq \psi(t_{0})-\psi(t)=\psi '(\xi)(t_{0}-t)=2\lambda \ip{f(\xi,\varphi(\xi)),\varphi(\xi)}\kr (t_{0}-t) <0\] for every $t\in [0,t_{0})$ and some $\xi\in(t,t_{0})$.
Thus, we assume that  $\psi(0)=R^{2}$. Furthermore, we estimate
\[ \Bigg|\jnt_{0}^{1}\varphi(s)\,dg(s)\Bigg|=\Bigg|\varphi(0)(g(0^{+})-g(0))+\lim_{\varepsilon\to 0^{+}}\jnt_{\varepsilon}^{1}\varphi(s)\,dg(s)\Bigg|>r_{-}.\]
Similarly, we obtain that 
\[ \Bigg|\jnt_{0}^{1}\varphi(s)\,dg(s)\Bigg| \leq r_{+}\] so Riemann-Stieltjes integral $\jnt_{0}^{1}\varphi(s)\,dg(s)$ satisfies the estimates in (ii). Then $h\Big(\jnt_{0}^{1} \varphi(s)\,dg(s)\Big)\neq 0$. Since $\varphi$ is the solution of (\ref{eq:FBVP}), we have a contradiction. \par 
According to the description of projections $P$ and $Q$ we have 
\[ (QN)(x)(t)=\Bigg(-h\Big(\jnt_{0}^{1}x(s)\,dg(s)\Big),h\Big(\jnt_{0}^{1}x(s)\,dg(s)\Big)\Bigg).\]
Since $\dim\ker{L}=\dim\im{Q}$, there exists an isomorphism $J:\im{Q}\to\ker{L}$. Let us define $J$ by
\[ J(-\pal,\pal)=\pal \quad (-\pal,\pal)\in \im{Q}.\] 
Then $(JQN)(x)=h\Big(\jnt_{0}^{1}x(s)\,dg(s)\Big)$. By the assumption (ii) we get that the topological Brouwer degree $\deg(JQN,\ker{L}\cap \Om,0)\neq 0$. Hence Mawhin theorem gives us the existence of the solution for (\ref{eq:ODE}), (\ref{eq:NC}) in the ball $\cl{B}_{C([0,1],\bb{R}^{k})}(0,R)$.  This completes the proof. \par
\textbf{Step 2.}\ Now, we assume that $\lim\limits_{\varepsilon\to 0^{+}}\var(g,[\varepsilon,1])\leq \min\limits_{j\leq k}|g^{j}(0^{+})-g^{j}(0)|$ and $\ip{f(t,x),x}\leq 0$ for $t\in(0,1]$, $|x|=R$ where $R>0$ is a constant. \par
 We consider the following BVP
\begin{equation}
\left\{\begin{array}{l}
x'=f(t,x)-\frac{1}{n}x, \\ 
h\Big(\jnt_{0}^{1}x(s)\,dg_{n}(s)\Big)=0, \quad n\in\bb{N},
\end{array}\right. \label{eq:PBVP}
\end{equation}
where $g_{n}=(g_{n}^{1},\ldots,g_{n}^{k}):[0,1]\to\bb{R}$ is such that $g_{n}(s)=g(s)$ for $s\in(0,1]$, $g^{j}_{n}(0)=g^{j}(0)$ for $j\neq j_{0}$ and $g_{n}^{j_{0}}(0)=g^{j_{0}}(0)-\frac{1}{n}\sgn(g^{j_{0}}(0^{+})-g^{j_{0}}(0))$, where $|g^{j_{0}}(0^{+})-g^{j_{0}}(0)|=\max\limits_{j\leq k}|g^{j}(0^{+})-g^{j}(0)|$. Then, functions $f(t,x)-\frac{1}{n}x$ and $g_{n}$ satisfy assumptions of Theorem \ref{thm:1thm}, so for every $n\in\bb{N}$ we get a solution of (\ref{eq:PBVP}) - $\varphi_{n}$. Moreover, $\norm{\varphi_{n}}_{C([0,1],\bb{R}^{k})}\leq R$ and sequence $(\varphi '_{n})_{n\in\bb{N}}$ is bounded in $C([0,1],\bb{R}^{k})$. Basing on the Ascoli-Arz\'ela theorem we can see that the sequence $(\varphi_{n})_{n\in\bb{N}}$ has a convergent subsequence in $C([0,1],\bb{R}^{k})$. We shall prove that the limit function $\varphi$ is solution of (\ref{eq:ODE}), (\ref{eq:NC}). Furthermore, since $\varphi_{n_{m}}$ is a solution of (\ref{eq:PBVP}), we have 
\[ \varphi '_{n_{m}}(t)=f(t,\varphi_{n_{m}}(t))-\frac{1}{n}\varphi_{n_{m}}(t)\to f(t,\varphi(t))\] uniformly as $m\to\infty$. Hence the limit function $\varphi$ is differentiable and $\varphi '(t)=f(t,\varphi(t))$. Let us observe that 
\[ \jnt_{0}^{1}\varphi^{j_{0}}_{n_{m}}(s)\,dg^{j_{0}}_{n_{m}}(s)=\varphi^{j_{0}}_{n_{m}}(0)(g^{j_{0}}_{n_{m}}(0^{+})-g^{j_{0}}_{n_{m}}(0))+\lim_{\varepsilon\to 0^{+}}\jnt_{\varepsilon}^{1}\varphi^{j_{0}}_{n_{m}}(s)\,dg^{j_{0}}_{n_{m}}(s)\]
\[=\varphi^{j_{0}}{n_{m}}(0)\big(g^{j_{0}}(0^{+})-g^{j_{0}}(0)\big)+\frac{1}{n_{m}}\sgn(g^{j_{0}}(0^{+})-g^{j_{0}}(0))\kr\varphi^{j_{0}}_{n_{m}}(0)+\lim_{\varepsilon\to 0^{+}}\jnt_{\varepsilon}^{1}\varphi^{j_{0}}_{n_{m}}(s)\,dg^{j_{0}}_{n_{m}}(s)\]
\[=\jnt_{0}^{1}\varphi^{j_{0}}_{n_{m}}(s)\,dg^{j_{0}}(s)+\frac{1}{n_{m}}\sgn(g^{j_{0}}(0^{+})-g^{j_{0}}(0))\kr\varphi^{j_{0}}_{n_{m}}(0).\]
Therefore $\jnt_{0}^{1}\varphi_{n_{m}}(s)\,dg_{n_{m}}(s)\to \jnt_{0}^{1}\varphi(s)\,dg(s)$ as $m\to\infty$. By the continuity $h:\bb{R}^{k}\to\bb{R}^{k}$ we obtain that $h\Big(\jnt_{0}^{1}\varphi(s)\,dg(s)\Big)=0$. Consequently $\varphi$ is a solution of (\ref{eq:ODE}), (\ref{eq:NC}).
\end{proof}  
\begin{rem} Let us assume that $g(1^{-})\neq g(1)$ and $\lim\limits_{\varepsilon\to 0^{+}}\var(g,[0,1-\varepsilon])\leq\min\limits_{j\leq k}|g^{j}(1)-g^{j}(1^{-})|$. By similar arguments the BVP (\ref{eq:ODE}), (\ref{eq:NC}) has at least one solution if there exists $\widehat{R}>0$ such that:
\begin{enumerate}
\item[(i')] $\ip{f(t,x),x}\geq 0$ for $t\in[0,1)$, $|x|=\widehat{R}$. 
\item[(ii')] Let \[\widehat{r}_{-}:=\widehat{R}\big(\min_{j\leq k}|g^{j}(1)-g^{j}(1^{-})|-\lim\limits_{\varepsilon\to 0^{+}}\var(g,[0,1-\varepsilon])\big),\] \[\widehat{r}_{+}:=\widehat{R}\big(|g(1)-g(1^{-})|+\lim\limits_{\varepsilon\to 0^{+}}\var(g,[0,1-\varepsilon])\big).\] $h(x)\neq 0$ for $\widehat{r}_{-}<|x|\leq \widehat{r}_{+}$ and the Brouwer degree $\deg(h,B_{\bb{R}^{k}}(0,\widehat{r}),0)$ is defined and does not vanish where $\widehat{r}\in(\widehat{r}_{-},\widehat{r}_{+}]$. 
\end{enumerate}
\end{rem}
\section{Applications}
Here we show the application of our results in the case of the second-order ordinary differential equation. We consider the following BVP
\begin{equation}
\left\{\begin{array}{l}
x''=f(t,x,x'), \\
h_{1}\Big(\jnt_{0}^{1}x(s)\,d\pog{g}_{1}(s),\jnt_{0}^{1}x'(s)\,d\pog{g}_{2}(s)\Big)=0, \\
h_{2}\Big(\jnt_{0}^{1}x(s)\,d\pog{g}_{1}(s),\jnt_{0}^{1}x'(s)\,d\pog{g}_{2}(s)\Big)=0,
\end{array}\right.\label{eq:BVPII}
\end{equation} 
where $f:[0,1]\times\bb{R}^{k}\times\bb{R}^{k}\to\bb{R}^{k}$, $h_{1}:\bb{R}^{k}\times\bb{R}^{k}\to\bb{R}^{k}$, $h_{2}:\bb{R}^{k}\times\bb{R}^{k}\to\bb{R}^{k}$ are continuous functions and $\pog{g}_{1}=(g
^{1},\ldots,g^{k}):[0,1]\to\bb{R}^{k}$, $\pog{g}_{2}=(g^{k+1},\ldots,g^{2k}):[0,1]\to\bb{R}^{k}$. It is obvious that problem (\ref{eq:BVPII}) is equivalent to BVP
\begin{equation}
\left\{\begin{array}{l}
\pog{x}'=\pog{f}(t,\pog{x}), \\
\pog{h}\Big(\jnt_{0}^{1}\pog{x}(s)\,dg(s)\Big)=0,
\end{array}\right.\label{eq:BVPI}
\end{equation}
where $\mathbf{x}=(x,y)$, $\pog{f}(t,\pog{x})=(y,f(t,x,y))$, $\pog{h}=(h_{1},h_{2})$, $y=x'$, $\pog{g}=(g^{1},\ldots,g^{k},g^{k+1},\ldots,g^{2k})$.
The problem (\ref{eq:BVPII}) has at least one solution if there exists $R>0$ such that \[\ip{x+f(t,x,y),y}\leq 0\] for $t\in(0,1]$, $|x|^{2}+|y|^{2}=R^{2}$, $h_{1}(x,y)\neq 0$, $h_{2}(x,y)\neq 0$ for $r^{2}_{-}<|x|^{2}+|y|^{2}\leq r^{2}_{+}$ where $r_{-}$, $r_{+}$ are defined in the theorem \ref{thm:1thm}, Brouwer degree $\deg\big((h_{1},h_{2}),B_{\bb{R}^{k}}(0,r)\times B_{\bb{R}^{k}}(0,r),0\big)$ is defined and does not vanish for some $r\in(r_{-},r_{+}]$ and a function $\pog{g}$ is an arbitrary function satisfying the assumptions of the theorem \ref{thm:1thm}. \par We will now discuss some special cases.  When $h_{1}$ depends only on $x$ and $h_{2}-y$, the condition of the function $h$ is as follows: degrees 
$\deg(h_{1},B_{\bb{R}^{k}}(0,r),0)$, $\deg(h_{2},B_{\bb{R}^{k}}(0,r),0)$ are defined and do not vanish. This is due to the following property (\cite{D}, page 33)
\[\deg((h_{1},h_{2}),B_{\bb{R}^{k}}(0,r)\times B_{\bb{R}^{k}}(0,r),0)=\deg(h_{1},B_{\bb{R}^{k}}(0,r),0)\kr \deg(h_{2},B_{\bb{R}^{k}}(0,r),0).\] \par
From now on we assume that $h_{1}(x,y)=x$ and $h_{2}(x,y)=y$. Moving away from the full generality, we assume that $\pog{g}_{1}=(g^{1},\ldots,g^{1})$, $\pog{g}_{2}=(g,\ldots,g)$, where 
\[ g^{1}(s)=\left\{\begin{array}{ccc}
1 & \text{for} & s=0, \\ 
0 & \text{for} & s\in(0,1],
\end{array}\right.\] and $g:[0,1]\to\bb{R}$ is an arbitrary function such that $\pog{g}=(\pog{g}_{1},\pog{g}_{2})$ satisfies the assumptions of theorem \ref{thm:1thm}. Hence we  obtain result for the problem (\cite{KT2})
\[ x''=f(t,x,x'), \quad x(0)=0, \quad \jnt_{0}^{1}x'(s)\,dg(s)=0,\] which is at resonance, then $g(1)-g(0)=0$. \par 
We now give another special case of BVP that we are generalized.
Namely, let us assume that $\pog{g}_{1}=(\widetilde{g},\dots,\widetilde{g})$, $\pog{g}_{2}=\pog{g}_{1}+\pog{g}_{3}$, where 
\[ \widetilde{g}(s)=\left\{\begin{array}{ccc}
1 & \text{for} & s=0, \\ 
0 & \text{for} & s\in(0,1],
\end{array}\right.\] and $\pog{g}_{3}=g=(g^{1},\ldots,g^{k}):[0,1]\to\bb{R}^{k}$ is an arbitrary function such that $\pog{g}=(\pog{g}_{1},\pog{g}_{2})$ satisfies the assumptions of theorem \ref{thm:1thm}. Hence we  obtain result for the problem 
\begin{equation} x''=f(t,x,x'), \quad x(0)=0, \quad x'(0)=\jnt_{0}^{1}x'(s)\,dg(s).\label{eq:P}\end{equation}
Similar problem was considered in \cite{SD}, with the difference that in second condition of (\ref{eq:P}) we have $x'(1)=\jnt_{0}^{1}x(s)\,dg(s)$. \par 
According to the introduction if $h_{1}(x,y)=x$, $h_{2}(x,y)=y$ and $\pog{g}_{1}=\pog{g}_{2}=(g,\ldots,g)$, where \[ g(s)=\left\{\begin{array}{ccc}
-1 & \text{for} & s=0,1, \\ 
0 & \text{for} & s\in(0,1),
\end{array}\right.\]
we obtain the result for classical periodic BVP (\cite{M}): $x(0)=x(1)$, $x'(0)=x'(1)$. However, it should be emphasized that our results do not embrace other classic BVPs such as the Dirichlet's problem and Neumann's problem.

\end{document}